\documentclass{aupl}

\usepackage{
amsfonts,
latexsym,
amssymb,
enumerate,
mathrsfs,
graphicx,
centernot,
color,
}

\newcommand{\labbel}[1]{\label{#1} [[{\bf #1}]]}  
\renewcommand{\labbel}{\label}

\usepackage{stackengine,xcolor,graphicx}
\stackMath
\def\dashsubseteq{\mathrel{
  \stackinset{l}{0pt}{c}{}{\wrule[2pt]{5pt}{.5pt}}{
  \stackinset{l}{1.6pt}{c}{}{\wrule[-10pt]{.5pt}{3pt}}{
  \stackinset{l}{0.5pt}{c}{1.5pt}{\rotatebox[origin=center]{45}{\wrule[2pt]{.5pt}{3pt}}}{
  \stackinset{l}{0.5pt}{c}{}{\rotatebox[origin=center]{-45}{\wrule[-5pt]{.5pt}{3pt}}}{
  \stackinset{l}{3.35pt}{c}{}{\wrule[-4pt]{.5pt}{10pt}}{
  \stackinset{l}{5.1pt}{c}{}{\wrule[-4pt]{.5pt}{10pt}}{
 \subseteq}}}}}}
}}
\newcommand\wrule[3][0pt]{\textcolor{white}{\rule[#1]{#2}{#3}}}

\newtheorem{theorem}{Theorem}[section]

\newtheorem{prop}[theorem]{Proposition} 
\newtheorem{proposition}[theorem]{Proposition} 
 
\newtheorem{corollary}[theorem]{Corollary}

\newtheorem*{claim*}{Claim}

\newtheorem*{theorem*}{Theorem}
\newtheorem*{proposition*}{Proposition}
\newtheorem*{corollary*}{Corollary}
\newtheorem*{lemma*}{Lemma}
\newtheorem*{scholion*}{Scholion}

\theoremstyle{definition}

\newtheorem{problem}[theorem]{Problem}

\theoremstyle{remark}
\newtheorem{remark}[theorem]{Remark}

\newtheorem*{remark*}{Remark}
\newtheorem*{remarks*}{Remarks}
 
\newtheorem{example}[theorem]{Example}

\newtheorem*{observation*}{Observation}

 \allowdisplaybreaks[1]

\numberwithin{equation}{section}

\begin{document}

\title[Comparable relations, amalgamation property]
{Comparable binary relations and the amalgamation property}

\author{Paolo Lipparini} 
\address{Diparti mento di Matematica\\Viale della  Ricerca
 Scientifica\\Universit\`a di Roma ``Tor Vergata'' 
\\I-00133 ROME ITALY}

\email{lipparin@axp.mat.uniroma2.it}

\subjclass{03C52; 06A75}

\keywords{Amalgamation property; Fra\"\i ss\'e limit,
comparable binary relations; relation-preserving operation}

\date{\today}

\thanks{
Work performed under the auspices of G.N.S.A.G.A. Work 
partially supported by PRIN 2012 ``Logica, Modelli e Insiemi''.
The author acknowledges the MIUR Department Project awarded to the
Department of Mathematics, University of Rome Tor Vergata, CUP
E83C18000100006.}

\begin{abstract}
 The theory of two binary
relations has the strong amalgamation property
when the first relation  is assumed to be coarser than
the second relation, and each relation satisfies 
a chosen set
of properties  from the following list:
transitivity, reflexivity, symmetry, 
antireflexivity and antisymmetry. 
 The amalgamation property is 
maintained when we add families of unary operations
preserving all the relations.
As a consequence, we get the existence of  Fra\"\i ss\'e 
limits for  classes of finite structures.

The results fail, for general comparability conditions,
 when three or more binary relations are
taken into account, or when we add an operation
preserving just one relation.
\end{abstract}

\maketitle

\section{Introduction} \labbel{intro} 

The amalgamation property has been proved useful
in algebra, model theory and algebraic logic, e.~g., 
\cite{CP,GG,H,KMPT,MMT}. 
Needless to say, binary relations are very important in
mathematics, henceforth it is interesting to study
which theories with one or more binary relations
have the amalgamation property.
Classical results furnish a positive answer for partial orders and
graphs \cite{F2}.

As a preliminary observation, we put the above results
in a unified context. We consider the following properties
of a binary relation: transitivity, reflexivity, symmetry, 
antireflexivity and antisymmetry, and 
we check that the class of structures with a binary 
relation satisfying any fixed subset of these properties has indeed
the amalgamation property.

More substantially, we show that the same holds
when we consider a comparable pair of binary relations,
even when the sets of properties satisfied by the two relations are distinct.
The amalgamation property is maintained 
when we add unary operations preserving all the relations.
The results are non trivial, since they fail, in general,
when we consider three or more relations,
or when we consider an operation preserving only one relation.

As a consequence, we get that the corresponding classes of finite structures
have a Fra\"\i ss\'e limit, and, in the absence of operations,
the corresponding theories have model completion.

\smallskip 

We now recall the main definitions. \emph{Structure}
is a synonym for first-order model \cite{H}. 
When we mention a class $\mathcal K$ of structures,
we shall always implicitly assume that all the structure
in the class have the same type (=language, =signature). For the sake of simplicity,
we shall deal with classes  closed under isomorphism.

 A class $\mathcal K$  closed under isomorphism
has the \emph{strong amalgamation property}, \emph{SAP},
 for short, if, whenever 
$\mathbf A, \mathbf B, \mathbf C \in \mathcal K$,
 $ \mathbf C \subseteq  \mathbf A$,
 $ \mathbf C \subseteq  \mathbf B$
and $C=A \cap B$,
then there is a structure
$\mathbf D \in \mathcal K$ such that 
$ \mathbf A \subseteq  \mathbf D$,
$ \mathbf B \subseteq  \mathbf D$.
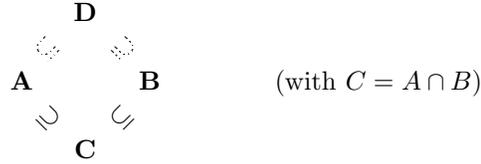
\begin{equation*}
\phantom{\qquad \qquad  \text{ (with $C=A \cap B$)}  }
 \begin{matrix} 
 \mathbf D \cr
$\rotatebox[origin=c]{45}{$ \dashsubseteq $}$
 \ \quad\  
$\reflectbox {\rotatebox[origin=c]{45}{$ \dashsubseteq $}}$ \cr
\mathbf A \quad \quad \quad \quad \mathbf B \cr
$\rotatebox[origin=c]{-45}{$ \supseteq $}$ \ \quad\ 
$\rotatebox[origin=c]{45}{$ \subseteq $}$
 \cr
 \mathbf C
 \end{matrix}   
\qquad \qquad  \text{ (with $C=A \cap B$)}  
\end{equation*}
We shall refer to $\mathbf A$, $\mathbf  B$, $\mathbf  C$ 
as a \emph{triple to be amalgamated}, a TBA-triple, for short. 
A theory $T$ has SAP if the class of models of $T$
has SAP, and similarly for the other properties we are going
to introduce.

A class has the \emph{amalgamation property}  (\emph{AP})
if the  conclusion in SAP is weakened to:
there are embeddings from $\mathbf A$ to $\mathbf  D$ 
and from $\mathbf B$ to $\mathbf  D$ which 
agree on $C$. In other words, under SAP, we can get $\mathbf  D$ 
in such a way that two elements $a \in A \setminus B$
and $b \in B \setminus A$ are never identified; this is not
necessarily the case if only AP is assumed.

If $R$ is a binary relation on some set, 
we frequently
write $ a \mathrel {R } b $
in place of $R(a,b)$ or $(a,b) \in R$.   
Moreover, $ a \mathrel {R } b \mathrel { S } c  $
is a shorthand for $ a \mathrel {R } b $ and $b  \mathrel { S } c  $.

In most cases, we shall prove a property somewhat stronger than 
SAP.
If $\mathcal K$ is a class of structures
closed under isomorphism and  in a language with a 
specified binary relation symbol $R$,
we say that $\mathcal K$ has the \emph{superamalgamation property},
or  \emph{superSAP},
\emph{with respect to $R$}\/
if, whenever 
$\mathbf A, \mathbf B, \mathbf C \in \mathcal K$
is a triple to be amalgamated, then 
there exists an amalgamating structure $\mathbf D
\in \mathcal K$  witnessing SAP  and such that,
for every  $a \in A \setminus B$ and $b \in B \setminus A$, 
  \begin{enumerate}[(i)]
    \item 
if $a \mathrel { R _{\mathbf D}} b $  
then there is $c \in C$ such that  $a \mathrel { R _{\mathbf A}} c 
\mathrel { R _{\mathbf B}} b $, and
\item
if $b \mathrel { R _{\mathbf D}} a $  
then there is $c \in C$ such that  $b \mathrel { R _{\mathbf B}} c 
\mathrel { R _{\mathbf A}} a $.
  \end{enumerate} 

We shall omit the reference to $R$ when the relation is understood.
The superamalgamation property
  is frequently  considered
with respect to  some
ordering relation \cite{GM}, but
there are  applications also in the case of an arbitrary binary relation
 \cite{Ma}.
The assumption that
$a \in A \setminus B$ and $b \in B \setminus A$
in the hypothesis of the superamalgamation properties
will simplify statements,
but notice that  many authors consider the weaker assumption 
$a \in A $ and $b \in B $.
Of course, the two definitions are equivalent if
$R$ is a reflexive relation, since if, say,
$a \in A \cap B$ and $a \mathrel { R} b$, then we can take $c=a$ in (i),
and similarly  for (ii).

A class $\mathcal K$ of structures  
for the same language has the \emph{joint embedding  property}
(JEP) if, for every  
$\mathbf A, \mathbf B  \in \mathcal K$,
 there are a structure
$\mathbf D \in \mathcal K$ and  embeddings
$ \iota: \mathbf A \rightarrowtail \mathbf D$
and 
 $ \kappa: \mathbf B \rightarrowtail \mathbf D$. 
In a sense, JEP can be considered as the special case of AP when $\mathbf  C$
is an empty structure, however, this does not always make sense,
for example 
if the language contains constants.

If  $\mathcal K$  is a class of finitely generated 
structures in a countable language, a  \emph{Fra\"\i ss\'e limit} 
of $\mathcal K$ is a 
countable ultrahomogeneous structure of age  $\mathcal K$.
 See \cite[Section 7.1]{H} for  details.
Classical examples are
 the ordered set of the rationals 
and the random graph which are, respectively,
 the Fra\"\i ss\'e
limit of the class of finite linearly ordered sets,
of the class of finite graphs.  

A 
first-order theory $T$  has \emph{model companion}
if $T$ has the same universal consequences of some model complete
theory $T^*$. If in addition $T$ has the amalgamation property,
the theory $T^*$
is a \emph{model completion} of $T$ \cite{H}.

\section{One relation with operators} \labbel{one}

Parts of the next proposition  follow from the 
proof of \cite[9.3]{F}. Compare also \cite[Lemma 3.3]{J}. 
Other cases are simple and well-known \cite{F2}.
The proposition provides a quite uniform treatment
for all the cases and deals also with additional operators,
a result which might be possibly  new.
If $\mathcal K$ is a class of structures,
$\mathcal K^{fin}$ 
denotes the class of the finite members of $\mathcal K$.

\begin{proposition} \labbel{pu}
Consider the following properties of 
a binary relation $R$.
\begin{enumerate}[1.]
    \item 
$R$ is transitive;
\item
$R$ is reflexive;
\item
$R$ is symmetric;
\item
$R$ is antireflexive, that is, 
$x \mathrel { R } x$ never holds;
\item
$R$ is antisymmetric.
 \end{enumerate}
Then the following statements hold.
  \begin{enumerate}[(A)]
    \item   
For every  $P \subseteq \{ 1,2,3,4,5\} $,
the class $\mathcal K_P$
 of the structures with a binary  relation
$R$ satisfying  the corresponding properties  has  superSAP.

\item
 For each $P \subseteq \{ 1,2,3,4,5 \} $, let $\mathcal K^f_P$
be the class of structures 
with an added unary operation $f$   
which is $R$-preserving, that is, 
\begin{equation}\labbel{rel}     
\text{ $  x \mathrel { R } y $ implies   
 $  f(x) \mathrel { R } f(y) $.}
 \end{equation}
Then  
$\mathcal K^f_P$
has superSAP.

\item
 For each $P \subseteq \{ 1,2,3,4,5 \} $, let $\mathcal K^g_P$
be the class of structures 
with an added unary operation $g$   
which is $R$-reversing, that is, 
\begin{equation}\labbel{rev}     
\text{ $  x \mathrel { R } y $ implies   
 $  g(y) \mathrel { R } g(x) $.}
 \end{equation}
Then  
$\mathcal K^g_P$
has superSAP.

\item
For every  $P \subseteq \{ 1,2,3,4,5 \} $ and
for every pair $F$, $ G$ of sets,
let $\mathcal K^{F,G}_P$
be the class of 
models obtained 
from members of $\mathcal K_P$
 by adding an  
 $F$-indexed set  of unary operations 
 satisfying  
\eqref{rel}
and a
$G$-indexed set of unary operations 
 satisfying  
\eqref{rev}.

Then
$\mathcal K^{F,G}_P$
has superSAP. 
 \end{enumerate}

For any one of the classes $\mathcal K$  considered in (A) - (D),
the class  $\mathcal K^{fin}$ 
has a Fra\"\i ss\'e limit, under the provision, in  (D), that the sets 
$F$ and $G$  are finite.
In case (A), for every  $P \subseteq \{ 1,2,3,4,5 \} $, 
if $\mathbf M$ is the Fra\"\i ss\'e limit of 
 $\mathcal K^{fin}_P$, then the
first-order theory $Th(\mathbf M)$ 
of $\mathbf M$ is
$ \omega$-categorical and has quantifier elimination; moreover, 
$Th(\mathbf M)$
is the 
model completion of $Th(\mathcal K_P)$.
  \end{proposition}

  \begin{proof}
Let us assume that 
$\mathbf A$, $\mathbf B$, $\mathbf C$
is a TBA triple (a Triple to Be Amalgamated).  
The proof of (A) will be divided into two  cases.

(A, case a) Suppose that  $ 1 \notin P$.
Then
let $R$ on $A \cup B$ be defined  
by $R= R _{\mathbf A} \cup R _{\mathbf B}$.
Namely, $d \mathrel { R } e $
if and only if 
\begin{equation}\labbel{ei}     
\text{ \underline{either} $d,e \in A$
and   $d \mathrel {  R _{\mathbf A}} e $,
\underline{or} 
 $d,e \in B$
and   $d \mathrel {  R _{\mathbf B}} e $.}
 \end{equation}
Let $\mathbf D = (A \cup B, R)$. 
Then, by construction,
$\mathbf A \subseteq \mathbf D$
and   
$\mathbf B \subseteq \mathbf D$, using the assumption that
$\mathbf C$
is a substructure of both $\mathbf A $ and $\mathbf B$.
If 
both $R_{\mathbf A} $ and $ R_{\mathbf B}$ 
are reflexive (symmetric, antireflexive, antisymmetric),
then so is $R$.
This is obvious for symmetry and  antireflexivity.
For reflexivity, use the assumption that
$D= A \cup B$.
 As far as 
antisymmetry is concerned, 
suppose that 
$R_{\mathbf A} $ and $R_{\mathbf B} $
are antisymmetric.
If $d \mathrel { R } e$, 
$e \mathrel { R } d$ and
these relations are  both witnessed by, say, $R_{\mathbf A} $, then 
 $d=e$, since  $R_{\mathbf A} $ 
is  antisymmetric.
On the other hand, if, say,  
$d \mathrel {  R _{\mathbf A}} e $
and   $e \mathrel {  R _{\mathbf B}} d $,
then necessarily $d,e \in A \cap B =C$, hence 
$d \mathrel {  R _{\mathbf C}} e $
and   $e \mathrel {  R _{\mathbf C}} d $,
since $\mathbf C$ is a substructure of both
$\mathbf A$ and $\mathbf B$.
Thus  $d=e$ in this case, as well, since  $R_{\mathbf C} $ 
is  antisymmetric.

In this case superamalgamation vacuously holds,
since if $a \in A \setminus B$ and $b \in B \setminus A$
then  $a$ and $b$ are not $R$-related. 

(A, case b)
Suppose that $ 1 \in P$, in particular, both
$R_{\mathbf A} $ and $R_{\mathbf B} $ are transitive.
Let $d \mathrel { R } e $
if  either \eqref{ei}, or
\begin{equation} \labbel{al}
\begin{aligned}
&d \in A, e \in B \text{ and there is $c \in C$ such that  }
d \mathrel { R_{\mathbf A} } c \mathrel { R_{\mathbf B} } e, \text{ or }  
\\  
&d \in B, e \in A \text{ and there is $c \in C$ such that  }
d \mathrel { R_{\mathbf B} } c \mathrel { R_{\mathbf A} } e.  
 \end{aligned}  
\end{equation}

Thus in this case we set
 $R= R _{\mathbf A} \cup R _{\mathbf B} \cup 
(R _{\mathbf A} \circ R _{\mathbf B}) \cup 
(R _{\mathbf B} \circ R _{\mathbf A})$.
Again, let $\mathbf D = (A \cup B, R)$. 

It is easy to verify that 
both $\mathbf A$ and $\mathbf B$ embed in  $\mathbf D$.
Say, if 
$d,e \in A$ and
 $d \mathrel { R } e$ is witnessed by
$d \mathrel { R_{\mathbf A} } c \mathrel { R_{\mathbf B} } e$,
for some $c \in C$, then necessarily
$e \in A \cap B =C$, thus 
 $c \mathrel { R_{\mathbf C} } e$, since 
 $\mathbf C$ embeds in  $\mathbf B$, hence
 $c \mathrel { R_{\mathbf A} } e$, since 
 $\mathbf C$ embeds in  $\mathbf A$, thus
$d \mathrel { R_{\mathbf A} } c \mathrel { R_{\mathbf A} } e$,
and
$d \mathrel { R_{\mathbf A} } e$
by transitivity of $  R_{\mathbf A}$.

It is not difficult to check that
$R$ is transitive.
In fact, we shall show that 
$R$ is the 
smallest transitive relation containing 
$R _{\mathbf A}$ and $R _{\mathbf B}$,
that is,
\begin{equation}\labbel{rrr}      
\text{$R= \bigcup _{n \geq 1}  R_1 \circ R_2 \circ  \dots \circ  R_n$,
\quad where each $R_k$ is either 
$R _{\mathbf A}$ or $R _{\mathbf B}$.}
  \end{equation}
Since $R _{\mathbf A}$ is transitive,
then $R _{\mathbf A} \circ R _{\mathbf A}
\subseteq R _{\mathbf A}$,
and similarly for $R _{\mathbf B}$,
hence we might assume that
$R _{\mathbf A}$ and $R _{\mathbf B}$
alternate in \eqref{rrr}.
Here we have used the fact that $\circ$
is associative and monotone
with respect to inclusion.  
If we also show that
$R _{\mathbf A} \circ R _{\mathbf B}\circ  R _{\mathbf A}
\subseteq R _{\mathbf A}$, then,
using also the symmetrical result, we get that in the union
in \eqref{rrr}  
it is enough to take $n \leq 2$,
whence the conclusion that $R$ is transitive. 
Indeed, let 
$ d \mathrel { R _{\mathbf A}} e 
 \mathrel { R _{\mathbf B}} g  \mathrel { R _{\mathbf A}} h$.
Then $e,g\in A \cap B = C$,
hence $e 
 \mathrel { R _{\mathbf B}} g $
implies 
$e 
 \mathrel { R _{\mathbf C}} g$,
since $\mathbf C \subseteq \mathbf B$.
But also   
 $\mathbf C \subseteq \mathbf A$, hence 
$e 
 \mathrel { R _{\mathbf A}} g $,
thus 
$d  \mathrel { R _{\mathbf A}} h$,
by transitivity of $R _{\mathbf A}$. 
We have proved that $R$ is transitive.

If 
both $R_{\mathbf A} $ and $ R_{\mathbf B}$ 
are reflexive (symmetric)
then trivially so is $R$.
For reflexivity, use again the assumption that
$D=A \cup B$. 
If $R_{\mathbf A} $ and $ R_{\mathbf B}$ 
are antireflexive, then so is  $R$,
since both $\mathbf A$ and $\mathbf  B$ embed
in $\mathbf  D$, and $D=A  \cup B$. 

We now deal with
 antisymmetry.
Suppose that $R_{\mathbf A} $ and $ R_{\mathbf B}$
are  antisymmetric,
$d \mathrel { R } e  $
and 
 $e \mathrel { R } d $,
the former relation witnessed, say, 
by $d \mathrel { R_{\mathbf A} } c \mathrel { R_{\mathbf B} } e$,
for some $c \in C$.
Since we have proved that $R$ is transitive,
then from $ c \mathrel { R_{\mathbf B} } e \mathrel { R } d $
we get $ c \mathrel { R } d $.
Since $c,d \in A$ and $\mathbf A \subseteq \mathbf D$, then 
$ c \mathrel { R_{\mathbf A} } d$, thus 
$c= d$, since   $R_{\mathbf A}$  
is antisymmetric.
Then from  $e \mathrel { R } d $ we get
 $e \mathrel { R } c $, hence 
$ e \mathrel { R_{\mathbf B} } c$, arguing as above, since
$e,c \in B$. But also    $ c \mathrel { R_{\mathbf B} } e $,
hence $e=c=d$, since  $  R_{\mathbf B}  $
is antisymmetric.
The other cases are similar or simpler.

In this case the superamalgamation
property follows immediately from the definition
of $R$.

We have concluded the proof of (A).

Given the proof of (A), clauses
(B), (C) and (D) are almost trivial.
If $\mathbf D$
extends both $\mathbf A$ and $\mathbf B$,
then  the operation $f$ is uniquely defined on
$D= A \cup B$ by
$f (a)= f _{\mathbf A}(a)$, for $a \in A$ and   
$f(a)= f_{\mathbf B}(a)$, for $a \in B$;
it is well-defined since $C= A \cap B$
and $\mathbf C$
is a substructure of both  
$\mathbf A$ and $\mathbf B$.
A similar comment applies to $g$.

If $R$ on $D$ is defined 
as in case a,
then \eqref{rel},
respectively, \eqref{rev}  obviously hold in $\mathbf D$.
If $R$ on $D$ is defined 
as in case b and 
$d \mathrel { R} e$
is witnessed, say, 
by 
$d \mathrel { R_{\mathbf A} } c \mathrel { R_{\mathbf B} } e$,
with $c \in C$, then 
$g_{\mathbf B}(e) \mathrel { R_{\mathbf B} } g_{\mathbf B}(c)=
g_{\mathbf A}(c) \mathrel { R_{\mathbf A} } g_{\mathbf A}(d)$,
since both $\mathbf A$ and $\mathbf B$
extend $\mathbf C$ and  satisfy \eqref{rev}
by assumption. 
Then  
$g(e) \mathrel { R_{\mathbf B} }g( c) \mathrel { R_{\mathbf A} } g(d)$,
hence $g(e) \mathrel { R } g(d)$, by \eqref{al}.  
The remaining cases are similar or simpler.

Notice that, as shown by the proof, it is not necessary
to assume that $F$ and $G$ 
are disjoint in (D).

To prove the last statements, first observe that,
for every class $\mathcal K$ under consideration,
the class $\mathcal K^{fin}$ has AP, since the proof
furnishes an amalgamating structure over $D=A \cup B$. JEP follows since
here we are allowed to consider an empty $\mathbf  C$.
Then use \cite[Theorems 7.1.2 and 7.4.1]{H}. 
 The assumption that $ F$ and $G$ are finite in (D) is necessary
in order to have a countable number of structures under
isomorphism.
 \end{proof} 

Neither the method applied in case a nor the method
applied in case b in the proof 
of Proposition \ref{pu} work for all situations.
See Remark \ref{div} below.
In fact, Proposition  \ref{div+}(a)(b) below implies that there is no
reasonable  ``monotone'' definition of $R$ 
on $\mathbf  D$   which covers all possible cases
in Proposition \ref{pu}.

\begin{remark} \labbel{str} 
Suppose that
$(A, \leq_{A})$,
$(B, \leq_{ B})$,
$(C, \leq_{C})$
is a TBA triple of partial orders
and $(D, \leq_{D})$
is obtained 
from the construction in case b in the proof of \ref{pu}.
If we let
\begin{equation}\labbel{r'}
e < f \quad \text{ if } \quad e \leq f \text{ and }
e \neq f  
   \end{equation}    
in each model, then 
$(D, <_{D})$
is obtained  applying case b to the models
$(A, <_{A})$,
$(B, <_{ B})$,
$(C, <_{C})$.

In words, the construction in \ref{pu}   
commutes in passing from some partial order
to the corresponding strict order.
In particular, we can equivalently deal with
orders in the strict or in the nonstrict sense.

Antisymmetry of $\leq$  is necessary in the above argument.
In fact, let $R$  be a  transitive  not antisymmetric
relation. If we define $R'$ by
\begin{equation*}\labbel{r''}
e \mathrel { R' } f \quad \text{ if } \quad e \mathrel { R }  f \text{ and }
e \neq f , 
   \end{equation*}     
then $R'$ might turn out to be not transitive, hence
case b cannot even be applied.  See
Proposition  \ref{rt} for a related counterexample. 
\end{remark}

If $P$ is a poset, a function $f:P \to P$
is \emph{strict order preserving (strict order reversing)}
if, respectively,
\begin{align*}
a < b \quad \text{ implies } \quad f(a) < f(b),
\\
a < b \quad \text{ implies } \quad  f(b) < f(a), 
  \end{align*}   
for all $a,b \in P$.

\begin{corollary} \labbel{propj}
Each of the following classes has superSAP and JEP.
Moreover, in each case, the class of finite structures
has a Fra\"\i ss\'e limit, provided the sets $F$ and $G$ 
below are finite.
 \begin{enumerate}
\item
The classes of 
  \begin{enumerate} [ \quad (a)]   
 \item 
 partially ordered sets,
\item  
preordered sets,
\item
undirected graphs (sets with a symmetric antireflexive binary relation),
\item
directed graphs (sets with an antireflexive binary relation),
\item
sets with  an equivalence relation,
\item
 sets with  a binary transitive relation,
\item
 sets with  a binary symmetric and reflexive  relation (a \emph{tolerance}),
  \end{enumerate}
 even when
an $F$-indexed family of relation-preserving \eqref{rel}
unary operations and 
a $G$-indexed family of  relation-reversing \eqref{rev}
unary operations are added.
\item
The class of partially ordered sets with
further families of order preserving, order reversing, strict order preserving
 and
strict order reversing unary operations.
 \end{enumerate} 
 \end{corollary}

\begin{proof}
Immediate from Proposition \ref{pu}.
Item (2) 
can be easily proved directly (e.~g., \cite[Corollary 2.4]{lop});
anyway, it follows from Remark \ref{str}
and Proposition \ref{pu} again.
 \end{proof}

A generalization of item (2) above in the strict case
to relations which are not orders
fails: see Proposition  \ref{rt}.

\section{Two  binary relations
with a comparability condition.} \labbel{ex2}

Given two binary relations $R$ and $S$ 
on the same domain,
we say that  $S$
is \emph{coarser} than  $R$ if $R \subseteq S$,
more explicitly, if $a \mathrel { R } b $ implies 
$a \mathrel { S } b $, for all $a$ and $b$ in the domain.
If this is the case, we shall also say that $R$ is \emph{finer}
than $S$. Notice that we shall always  use the expression
``coarser'' in the sense of ``coarser than or equal to''.

Recall the properties 1.- 5. from 
Proposition \ref{pu}. 
They correspond, in this order, to 
transitivity, reflexivity, symmetry, 
antireflexivity and antisymmetry.

\begin{theorem} \labbel{due} \ 
  \begin{enumerate}[(A)]
    \item   For every pair  $P, Q \subseteq \{ 1,2,3,4,5 \} $,
the class $\mathcal K_{P,Q}$ of structures with 
  \begin{enumerate}
   \item  
a binary relation $R$ satisfying the properties 
from $P$ and
\item 
 a coarser
relation $S$ satisfying the properties 
from $Q$
  \end{enumerate}
 has SAP.
Actually, superSAP
holds simultaneously both with respect to 
$R$ and $S$.

The subclass of finite structures of
$\mathcal K_{P,Q}$ has a Fra\"\i ss\'e limit 
 $\mathbf M$. 
The
first-order theory $Th(\mathbf M)$ is
$ \omega$-categorical and has quantifier elimination; moreover, 
$Th(\mathbf M)$
is the 
model completion of the
 theory $Th(\mathcal K_{P,Q})$.

\item 
 SuperSAP is maintained if we add families of 
  \begin{enumerate}[(i)]
    \item 
unary operations which are both $R$- and $S$-preserving;
\item
unary operations which are both $R$- and $S$-reversing;
\item
unary operations which are  $R$-preserving;
\item
unary operations which are $R$-reversing;
  \end{enumerate} 

Fra\"\i ss\'e limits of finite structures
still exist, provided only a finite number of operations are added.

\item
 On the other hand, 
the class of structures with a transitive relation $R$,
a coarser binary relation $S$ and an $S$-preserving function $f$  
has not AP (here we do not include the condition that $f$ 
is  $R$-preserving).

\item
Similarly, the class of structures with a partial order $\leq$,
a coarser symmetric and reflexive relation $S_1$ and an 
$S_1$-preserving function $f$  
has not AP.
\end{enumerate} 
 \end{theorem}

  \begin{proof} 
 (A)
The result follows in an obvious way from
the proof of Proposition \ref{pu}
when either case a or case b can be applied to both
$R$ and $S$, since \eqref{ei} and \eqref{al}
preserve coarseness.
Similarly, the result holds when
case a is applied to $R$ and case b is applied to $S$,
since case b always provides a coarser relation
in comparison with  case a.

Hence we can suppose that $1 \in P$ and 
$1 \notin Q$, that is, $R$ is transitive and
$S$ is not supposed to be transitive.
If $\mathbf A$, $\mathbf  B$, $\mathbf  C$ 
is a TBA triple, let $D=A \cup B$
and define $R$ on $D$ as in case b
in the proof of Proposition \ref{pu},
namely, $ d \mathrel { R} e $
holds if either \eqref{ei} or \eqref{al} applies.  
By the proof of Proposition \ref{pu},
$R$ on $D$ inherits from $\mathbf A$ and 
$\mathbf  B$  all the properties in $P$. 
 
It remains to define $S$ on $D$.
Two cases need to be considered.
Case $3 \in Q$. This means that $S$ is required to be symmetric. 
In this case we define $S$ on $D$ by  $ d \mathrel { S} e $ if 
\begin{equation}\labbel{als}   
\text{ \underline{either} $d,e \in A$
and   $d \mathrel {  S _{\mathbf A}} e $,
\underline{or} 
 $d,e \in B$
and   $d \mathrel {  S _{\mathbf B}} e $,
\underline{or}
$d \mathrel {  R } e $, 
\underline{or}
$e \mathrel {  R } d $.}
  \end{equation} 

Observe that  $S$ is symmetric and coarser than $R$ by construction.
Then notice that if $d, e \in A$
and $d \mathrel {  R } e $, 
then $d \mathrel {  R _{\mathbf A}} e $,
since, as proved in Proposition \ref{pu},
$(A,  R _{\mathbf A})$ embeds
in $(D,R)$, hence also 
$d \mathrel {  S _{\mathbf A}} e $,
since 
$R$ is finer than $S$ on $\mathbf A$.
A similar remark applies if 
$d, e \in B$.
Thus, when applying the two last disjuncts
in \eqref{als},
we can restrict ourselves to the case 
when $d\in A \setminus B$
and $e \in B \setminus A$,
or conversely. 
This observation immediately implies that both
$(A,  S _{\mathbf A})$ and
$(B,  S _{\mathbf B})$
 do embed in $(D,S)$,
in particular, if $S$ is reflexive
(antireflexive) on $\mathbf A$ and $\mathbf  B$,
then  $S$ is reflexive
(antireflexive) on $D$.

Expanding the argument a bit, we get  superSAP with respect to $S$ (notice that superSAP with respect to $R$ follows from Proposition \ref{pu}).
Indeed, if, say, $d\in A \setminus B$,
 $e \in B \setminus A$ and 
$d \mathrel {  S } e $, then
either $d \mathrel {  R } e $
or $e \mathrel {  R } d $, by the definition of $S$.
Suppose, say, the latter holds.
Then, by \eqref{al},
there is $c \in C$ such that 
 $e \mathrel { R_{\mathbf B} } c \mathrel { R_{\mathbf A} } d$.
Since $R$ is finer than $S$, we get
 $e \mathrel { S_{\mathbf B} } c \mathrel { S_{\mathbf A} } d$,
hence  $d \mathrel { S_{\mathbf A} } c \mathrel { S_{\mathbf B} } e$,
by symmetry of $S_{\mathbf A} $ and $  S_{\mathbf B}$. 
Then the definition of $S$ provides
 $d \mathrel S c \mathrel S e$, 
what we had to show.

The case when $5 \in Q$, that is,
 $S$ is required to be antisymmetric, is trivial, since
no two distinct elements are $S$-connected by
a symmetric and antisymmetric relation
(Recall that we have already proved that
  both $\mathbf A$  and $\mathbf  B$  embed in $\mathbf  D$,
 and $D = A \cup B$.)
We have showed that if 
$3 \in Q$, then all the properties 
in $Q$ are preserved.

Case $3 \notin Q$. This means that $S$ is not required to be symmetric.
In this case we define $S$ on $D$ by  $ d \mathrel { S} e $ if 
\begin{equation}\labbel{alss}   
\text{ \underline{either} $d,e \in A$
and   $d \mathrel {  S _{\mathbf A}} e $,
\underline{or} 
 $d,e \in B$
and   $d \mathrel {  S _{\mathbf B}} e $,
\underline{or}
$d \mathrel {  R } e $.}
  \end{equation} 

As in case $3 \in Q$, 
$S$ is  coarser than $R$, and 
\eqref{alss} provides $S$-embeddings.
The only nontrivial case left to check is when 
$S$ is required to be antisymmetric
and the only nontrivial part is when 
$d \mathrel {  S } e $,
$e \mathrel {  S } d $
and, say, $d \in A \setminus B$
and $e \in B \setminus A$.
Thus 
$d \mathrel {  S } e $
and $e \mathrel {  S } d $
are witnessed by
$d \mathrel {  R } e $,
$e \mathrel {  R} d $,
hence, since $R$ on $D$ is defined by 
\eqref{al}, there are $c, c' \in C$
such that   
$d \mathrel { R_{\mathbf A} } c \mathrel { R_{\mathbf B} } e$
and
$e \mathrel { R_{\mathbf B} } c' \mathrel { R_{\mathbf A} } d$.
Then
$c \mathrel { R_{\mathbf B} } e \mathrel { R_{\mathbf B} } c'$,
hence 
 $c \mathrel { R_{\mathbf B} }  c'$, since
$ R_{\mathbf B} $
is transitive, thus
 $c \mathrel { S_{\mathbf B} }  c'$,
since 
$ \mathrel { R_{\mathbf B} }  $ 
is finer than
$ \mathrel { S_{\mathbf B} }  $;
also 
$c \mathrel { S_{\mathbf C} }  c'$,
since $c, c' \in C$ and $\mathbf  C$ embeds in $\mathbf  B$.
Similarly, $c' \mathrel { S_{\mathbf C} }  c$,
hence $c=c'$, since  $ \mathrel { S_{\mathbf C} }$ 
is antisymmetric.
Hence $c \mathrel { R_{\mathbf B} } e \mathrel { R_{\mathbf B} } c$,
thus 
$c \mathrel { S_{\mathbf B} } e \mathrel { S_{\mathbf B} } c$, hence
$c=e$, by antisymmetry of $S_{\mathbf B}$;
similarly, $c=d$, hence $e=d$, what we had to show.    

The statement about Fra\"\i ss\'e limits
and model completions is proved as 
Proposition \ref{pu}, again using  \cite[Theorems 7.1.2 and 7.4.1]{H}. 

In order to prove (B) (i), if $f$ is an $R$- and $S$-preserving
unary function on a TBA triple   $\mathbf A$, $\mathbf  B$, $\mathbf  C$,
define $f$ in the unique compatible way on 
$D= A \cup B$. By the proof of Proposition \ref{pu}, $f$
is $R$-preserving; moreover, the proof of \ref{pu}
entails $S$-preservation, too, 
apart from the exceptional case  $1 \in P$ and 
$1 \notin Q$. It remains to show that 
$f$ is $S$-preserving in this case, too. This is obvious if 
$ d \mathrel { S} e $ is given by    
  $d \mathrel {  S _{\mathbf A}} e $
or 
   $d \mathrel {  S _{\mathbf B}} e $
 in \eqref{als} or in  \eqref{alss}.     
If 
$ d \mathrel { S} e $ is given by, say,    
$ e \mathrel { R} d $ in \eqref{als},
then we know that 
$ f(e) \mathrel { R} f(d) $.
But then 
$ f(d) \mathrel { S} f(e) $,
applying 
 \eqref{als}. The other
cases are 
slightly 
simpler.

 (B) (ii) is entirely similar.

(B) (iii) and (iv) are immediate from the proof of
Proposition  \ref{pu}(D).

(C) Let $C= \{  c_1, c_2, c_3\} $
with no pair of elements $R$-related, 
$ c_1 \mathrel { S }  c_2   \mathrel { S }   c_3$,
no other $S$-relation and $f$ the identity function.
Extend $\mathbf  C$ to $\mathbf A$   
 by letting $A= C \cup \{  a\} $,
$a \notin C$, $a \mathrel { R  } c_2$,
$a \mathrel { S  } c_2$ and no other 
$R$- or  $S$-relation, apart from the relations coming from
the assumption that $\mathbf  C$ embeds in $\mathbf A$.
Let $f(a)=c_1$. 
 Similarly,  let $B= C \cup \{  b\} $
with $A \cap B =C$, 
 $c_2 \mathrel { R  } b$,
$c_2 \mathrel { S  } b$, no other 
unnecessary $R$- or  $S$-relation
and $f(b)=c_3$.

In any amalgamating structure
we should have 
 $a \mathrel { R  } b$ by transitivity of $R$,
hence $a \mathrel { S  } b$,
by the coarseness assumption.
But then $S$-preservation of  
$f$ implies $c_1= f(a) \mathrel { S } f(b)=c_3$,
a contradiction. 

A counterexample for 
(D) is obtained in a similar way, e.~g.,
let $\leq$ be the smallest partial order containing $R$ 
as defined in (C) above, and take $S_1$ here to be the smallest 
reflexive and symmetric relation containing $S$ 
as defined in (C).
\end{proof} 

The analogue of Theorem  \ref{due}(A)
when three or more relations are taken into account is false.
See Proposition \ref{div+}(a)(b).
See Proposition \ref{div+}(c)  for a counterexample
related to \ref{due}(C).

In particular, it follows from Theorem \ref{due}
that if $R$ is a binary relation supposed
to satisfy some condition from clauses (a)-(g)
in Corollary \ref{propj}(1) and $S$ is another relation 
supposed to satisfy a possibly distinct condition from (a)-(g), then 
the class of all structures in which $S$ is coarser than
$R$ has SAP. 
We single out a few special cases of interest.

\begin{corollary} \labbel{cordue}   
The following theories have
superSAP with respect to  both the relations
involved.

  \begin{enumerate}
    \item  
The theory of  two
partial orders, one coarser than the other.
    \item  
The theory of a 
partial order with a coarser preorder relation.
\item
The theory  of a partial order with a coarser tolerance 
 relation.
\item
The theory  of a strict partial order with a coarser graph 
 relation.
\item
The theory of a directed graph
with a coarser strict order relation.
  \end{enumerate}

SAP is maintained if families of both-relations preserving or reversing
unary operations are added.
\end{corollary}

\begin{remark} \labbel{rimpo}
The theories in Corollary \ref{cordue}(1)(2) 
 have found many applications in very disparate
fields, such as algebraic geometry, domain theory,
information systems,  foundations of topology and of
general relativity. See \cite{mtt} for more details and references. 
\smallskip 

The theory $T$  considered in Corollary \ref{cordue}(5)
seems to be interesting.
Clearly,  the graph reduct of any model  of $T$  
is acyclic. On the other hand, the class of acyclic
graphs has not AP, though every acyclic graph
can be expanded (not in a unique way) to a model of $T$. 
The contradiction is only apparent, we explain the point 
by an example.

Let $\mathbf  C$ be the directed graph 
with $C= \{  c, d\} $ and the elements not connected.
Let $A= C \cup \{ a \} $
with an edge from $c$ to $a$   
and an edge from $a$ to $d$.
Let $B= C \cup \{ b \} $
with an edge from $d$ to $b$   
and an edge from $b$ to $c$.
Clearly, $\mathbf A$ and $\mathbf  B$ cannot be
amalgamated over $\mathbf  C$ in the class
of acyclic directed graphs.

Now we try to endow $\mathbf A$, $\mathbf  B$
and $\mathbf  C$ with a coarser order relation.
We have no constraint on $\mathbf  C$.
However, by the coarseness assumption,
we must have $c<a<d$ in the expanded $\mathbf A$.  
This does not furnish a contradiction: 
if $c<a<d$, hence $c<d$  in $\mathbf A$,
we must have $c<d$ in (a possible expansion of) $\mathbf C$,
since we want  $\mathbf  C$ to be a substructure of $\mathbf A$.

Arguing symmetrically,  $d<b<c$ in $\mathbf  B$
so if we want $\mathbf  C$ to be embedded in  $\mathbf B$,
we must have $d < c$. 
Each of the above conditions, taken alone, can be realized,
but we cannot have both conditions holding at the same time.
Thus there is no way
of adding an order relation to $\mathbf  C$  
 to make it a substructure of both the expanded $\mathbf A$ and
$\mathbf  B$ as models of $T$.

What have we got?
We can turn the class of acyclic graphs into a class
with SAP, provided we add a further order relation
which puts some constraints on pairs of elements allowed
to be connected ``in the future''  by some path in some extension. 
 
In the above respect, we remark that 
expanding the language is a well-known strategy in model
theory. We do not know whether the idea has already been used 
before in 
the present situation.
\end{remark}   

\begin{problem} \labbel{apexp}  
The possibility of turning some theory without AP
into a theory with AP by means of expansions, as exemplified
in Remark \ref{rimpo}(b),  seems rather interesting.

Are there applications of the theory 
in Corollary \ref{cordue}(5)
and of the remarks in  \ref{rimpo}?
\end{problem}

We now show that Theorem \ref{due} 
does not generalize to the case of three relations.

\begin{prop} \labbel{div+}
The following theories do not have  AP.
  \begin{enumerate}[(a)]
    \item   
The theory of an antisymmetric relation $S$ 
with two partial orders $\leq$ and $\leq'$
both finer than $S$. 

\item
 The theory of an antisymmetric relation $S$ 
with two transitive relations
both finer than $S$.

\item
 The theory of a partial order with
a coarser antisymmetric relation $S$
and a bijective $S$-preserving 
unary operation.  
\end{enumerate} 
 \end{prop}   

\begin{proof} 
(a) Let $C= \{  c,d\} $
with no pair $S$-related and with 
$c$ and  $d$ both   $\leq$- and $\leq'$-incomparable.
Let $A= C \cup \{  a \} $ with 
$a \notin C$, $c \leq a $, $a \leq' d$,     
  $c \mathrel { S }  a $, $a \mathrel { S}  d$
and no other $S$- or order-relation, apart
from order reflexivity.
Let $B= C \cup \{  b \} $ with 
$A \cap B = C$, $b \leq c $, $d \leq' b$,     
  $b \mathrel { S }  c $, $d \mathrel { S}  b$
and no other order- or $S$-relation, apart
from order reflexivity.
In any amalgamating structure satisfying the theory we must have
$b \leq a$, hence  $b \mathrel { S }  a $,
since $\leq$ is finer thn $S$.
Similarly, $a \leq ' b$, hence 
$a \mathrel { S }  b$, but this contradicts 
antisymmetry of $S$.
We cannot have $a=b$ since, say,
 $a \mathrel { S}  d$ but not $b \mathrel { S}  d$ and both
$\mathbf A$ and $\mathbf  B$ should be embedded in $\mathbf  D$.

The proof of (b) presents no significant difference.

(c) Let $C= \{ c_1,d_1, c_2,d_2\} $
with all the elements  pairwise $\leq$-incomparable
 and $S$-unrelated. Here and below let $f$ be the bijection
which permutes the indices.

Extend $\mathbf  C$ to $\mathbf A$ with
$A= C \cup \{ a_1, a_2 \}$,
 $a_1 \geq c_1$, $d_2 \geq a_2$   
and 
 $d_i \mathrel { S }  a_i \mathrel { S }  c_i$
($i=1,2$). 
Extend $\mathbf  C$ to $\mathbf B$ with
$B= C \cup \{ b_1, b_2 \}$,
 $c_1 \geq b_1$, $b_2 \geq c_2$  
and 
 $c_i \mathrel { S }  b_i \mathrel { S }  d_i$
($i=1,2$). 

In any amalgamating structure
we should have $a_1 \geq b_1$ and 
 $b_2 \geq a_2$, by transitivity of $\geq$.
Then $a_1 \mathrel { S }  b_1$ and 
 $b_2 \mathrel { S } a_2$, by the coarseness assumption.
Since $f$ is  $S$-preserving, then
 $b_1=f(b_2) \mathrel { S }  f(a_2)=a_1$,
hence $a_1=b_1$ since $S$ is supposed to be
antisymmetric. 
But then $c_1=a_1=b_1$,
since $\leq$ should be a partial order,
a contradiction. 
 \end{proof}

It is possible to  show 
that Theorem \ref{due} 
 generalizes to any number of relations, 
provided all the relations under consideration
are transitive.
Quite unexpectedly, a  generalization
of  Theorem \ref{due} 
does hold 
 when dealing with many relations, some of them
possibly non transitive.
We shall show that the only obstacle to 
amalgamation are antisymmetric relations
supposed to be  coarser than a pair of incomparable transitive relations;
essentially, the only counterexamples to amalgamation
are given by the examples we have described
 in Proposition \ref{div+}.
However,
in case we add operations, the operations are required to preserve all the relations,
compare Theorem \ref{due}(C)(D)
and Proposition \ref{div+}(c).  
We shall present  details elsewhere.

\section{Further remarks} \labbel{fur}

\begin{remark} \labbel{div}
Neither the method applied for case a nor the method
applied for case b in the proof 
of Proposition \ref{pu} work for all situations.

In general, if transitivity is not assumed, it might happen
 that $c \mathrel {R } d $,
 $d \mathrel {R } e$
but not $c \mathrel {R } e $
in $\mathbf C$.
Then, no matter the structure of 
$\mathbf A$ and $\mathbf B$,
if we apply clause \eqref{al},
we get  $c \mathrel {R } e $ in $\mathbf D$,
hence $\mathbf C$ is not even a substructure of 
$\mathbf D$.  

Even worse! It might happen that 
$\mathbf A$ and $\mathbf B$
are actually substructures of 
$\mathbf D$, but antisymmetry is not preserved.
Let $C= \{ c,d \} $, with no pair of elements 
$R$-related in $\mathbf C$.
Let $c \mathrel {R } a  \mathrel {R } d $
in $\mathbf A$ and 
$d \mathrel {R } b  \mathrel {R } c $
in $\mathbf B$.
Then $a  \mathrel {R _{\mathbf A}} d \mathrel {R _{\mathbf B}} b$,
thus if clause  \eqref{al} holds in $\mathbf D$,
we have $a \mathrel { R} b $ and symmetrically
$b \mathrel { R} a $.
Since we are asking SAP,
$a \neq b$, hence antisymmetry fails in 
  $\mathbf D$.

On the other hand, concerning case b, if we want transitivity to be preserved,
clause \eqref{ei} alone is clearly not sufficient. 
 \end{remark}   

We now see that not all possible variations on 
Proposition \ref{pu} hold. 

\begin{proposition} \labbel{rt}
Let $T$ be the theory  
with a binary reflexive and transitive relation $R$
and a unary function $f$  which strictly preserves $R$,
namely, 
\begin{equation*}    
\text{$d \mathrel { R } e $
and $d \neq e$ \quad imply \quad
both  
$f(d) \mathrel { R } f(e) $
and $ f(d) \neq f(e)$.}
 \end{equation*} 

Then $T$ has not AP. 
\end{proposition}   

\begin{proof}
Let $C= \{  c, d,e  \}$
with three distinct elements such that $ d \mathrel { R } e \mathrel { R } d $,
$ x \mathrel { R } x$, for all $x$, 
and    $f(d)=e$, $f(e)=d$ and $f(c)=d$.   
Let $\mathbf A$ extend $\mathbf C$
with   
$A= \{ a, c, d,e  \}$, $a \notin C$,
 $ a \mathrel { R } c $, $a \mathrel { R } a $
and    $f(a)=e$.
Let $\mathbf B$ extend $\mathbf C$
with   
$B= \{ b, c, d,e  \}$,
$b \notin C$, $b \neq a$, 
 $ c \mathrel { R } b $, $b \mathrel { R } b$
and    $f(b)=e$.

If $\mathbf D$ strongly amalgamates 
 $\mathbf A$ and $\mathbf B$ over $\mathbf C$,
then $ a \mathrel { R } b $
 by transitivity, but then
$f(a)=e=f(b)$,
contradicting  
strict preservation of $R$.
Hence SAP fails.
We cannot identify $a $ and $ b$ in $\mathbf D$
(hence we cannot even have AP),
  since $ c \mathrel { R } b $
in $\mathbf B$, but 
 $ c \mathrel { R } a $ is not assumed to hold in $\mathbf A$
and such relationships must be preserved in any common embedding.
\end{proof}
 
Notice that, on the other hand, the theory $T$  
with a binary reflexive and transitive relation $R$
and a \emph{bijective} function $f$  which  preserves $R$
has superSAP. This is immediate from the proof of 
Proposition \ref{pu}. 
Obviously, $R$-preserving bijective functions
are also strict  $R$-preserving.

\begin{example} \labbel{strpre}   
The theory 
with an equivalence relation $R$
and a unary function $f$  which strictly preserves $R$
has SAP but not AP. A counterexample to SAP is easy, just let
$a \in A \setminus C$ and $b \in B \setminus C$ be such that
$f(a)=f(b) \in C$. In order to prove AP,
given a TBA triple
$\mathbf A$, $\mathbf  B$, $\mathbf  C$,
identify $a \in A$ and $b \in B$ if
$a \mathrel { R  } c \mathrel { R } b$,
for some $c \in C$ and, moreover,
$f^n(a)=f^n(b) \in C$.  
 With the above identifications, 
$A \cup B$ becomes an appropriate amalgamating structure,
with $R $ the transitive closure of $R _{ \mathbf A} \cup R _{ \mathbf B}$
and $f$ defined in the unique compatible way.
Notice that if $a$ and $b$ are identified as above, then
$f(a)$ and  $f(b)$ are identified, as well, since 
either $n=1$ and $f(a)=f(b) $, or $n > 1$
and  $f^{n-1}(f(a))= f^n(a)=f^n(b) =f^{n-1}(f(b)) \in C $. 
Moreover, the assumption $a \mathrel { R  } c \mathrel { R } b$
assures that we do not identify distinct $R$-classes. 
 \end{example}

\section{Statements and declaration} \labbel{st} 

Ethical statement. This material is the author's own original work.
 The paper reflects the author's own research and analysis in a truthful and complete manner.
 All sources used are properly disclosed. This study has not involved human subjects.
This research did not involve the use of animals. 

The author has no competing interests to declare that are relevant to the content of this article.


\begin{thebibliography}{}    


\bibitem[CP]{CP} Czelakowski, J.,   Pigozzi, D.,
\emph{Amalgamation and interpolation in abstract algebraic logic},
 in Caicedo, X., Montenegro, C. H. (eds.), 
\emph{Models, algebras, and proofs ({B}ogot\'{a}, 1995)},
Lecture Notes in Pure and Appl. Math.,
\textbf{203},
187--265 (1999).





\bibitem[F1]{F}
Fra\"{\i}ss\'{e}, R.,
\emph{Sur l'extension aux relations de quelques propri\'{e}t\'{e}s des
ordres},
Ann. Sci. Ecole Norm. Sup. (3),
\textbf{71}, 363--388
(1954).

\bibitem[F2]{F2}
 Fra\"{\i}ss\'{e}, R.,
\emph{Theory of relations. Revised edition with an appendix 
by N. Sauer},
Studies in Logic and the Foundations of Mathematics \textbf{145},
 Amsterdam: North-Holland (2000). 




\bibitem[GM]{GM}
 Gabbay, D. M., Maksimova, L., 
\emph{Interpolation and definability. Modal and intuitionistic logics},
Oxford Logic Guides
\textbf{46}, The Clarendon Press, Oxford University Press, Oxford
(2005).

\bibitem[GG]{GG} Ghilardi, S.,  Gianola, A., 
\emph{Modularity results for interpolation, amalgamation and superamalgamation}, 
Ann. Pure Appl. Logic \textbf{169}, 731--754 (2018).




\bibitem[H]{H} 
Hodges, W., \emph{Model theory},
 Encyclopedia of Mathematics and its Applications \textbf{42},
 Cambridge University Press, Cambridge, 1993. 



\bibitem[J]{J} J{\'o}nsson, B.,
\emph{Universal relational systems},
Math. Scand. \textbf{4}, 193--208  (1956).   



\bibitem[KMPT]{KMPT} 
 Kiss, E. W., M\'arki, L., Pr\"ohle,  P.,  Tholen,  W.,
\emph{Categorical algebraic properties. 
A compendium on amalgamation, congruence extension, 
epimorphisms, residual smallness, and injectivity},
 Studia Sci. Math. Hungar. \textbf{18},  79--140  (1982). 


\bibitem[L]{lop}  Lipparini, P.,  
\emph{Linearly ordered sets with only one operator
have the amalgamation property},
Ann. Pure Appl. Logic \textbf{172},
Paper No. 103015, 15 (2021).

 \bibitem[L2]{mtt}  Lipparini, P.,  \emph{A model theory of topology}, arXiv:2201.00335, 1--30 (2022).




\bibitem[M]{Ma}  
 Madar\'{a}sz,  J. X.,
 \emph{Interpolation and amalgamation; pushing the limits. {I}},
Studia Logica
\textbf{61},
    311--345 (1998).



\bibitem[MMT]{MMT}  Metcalfe, G.,  Montagna,  F.,  Tsinakis, C., 
\emph{Amalgamation and interpolation in ordered algebras}, J. Algebra 
\textbf{402},  21--82 (2014).



\end{thebibliography}
\end{document}